\documentclass[12pt]{article}
 \usepackage{amsmath,amssymb}

\newcommand{\Var}{{\cal{V}_{\mathbb{C}}}}
\newcommand{\VarS}{{\cal{V}_{S}}}
\newcommand{\hilb}{{\mbox{\rm Hilb}}}

\def\uu{{\underline{u}}}

\def\tt{{\underline{t}}}
\def\ii{\underline{i}}

\def\nn{\underline{n}}

\def\kk{{\underline{k}}}

\def\1{\underline{1}}
\def\AA{{\mathbb A}}

\def\LLL{{\mathbb L}}
\def\Z{{\mathbb Z}}
\def\C{{\mathbb C}}

\def\ZZ{{\mathcal Z}}

\def\CP{\mathbb C\mathbb P}

\def\gm{{\mathfrak m}}
\def\gM{{\mathfrak M}}
\def\glocal{{\mathfrak M}_{\C^{d-1}\times \C,\{O\}\times \C,O}}

\newtheorem{theorem}{Theorem}

\newtheorem{proposition}{Proposition}

\newenvironment{definition}
{\smallskip\noindent{\bf Definition\/}:}{\smallskip\par}

\newenvironment{example}
{\smallskip\noindent{\bf Example\/}.}{\smallskip\par}

\newenvironment{remark}
{\smallskip\noindent{\bf Remark\/}.}{\smallskip\par}

\newenvironment{proof}
{\noindent{\bf Proof\/}.}{{ $\square$}\smallskip\par}

\title{On the power structure over the Grothendieck ring of varieties and its applications
\footnote{Math. Subject Class.: 14C05, 14G10}
}
\author{S.M.~Gusein-Zade \thanks{Partially supported by the grants
RFBR--04--01--00762,
NSh--4719.2006.1, and NWO-RFBR 047.011.2004.026.
Address: Moscow State University,
Faculty of Mathematics and Mechanics, Moscow, 119992, Russia.
E-mail: sabir\symbol{'100}mccme.ru} \and
I.~Luengo
\and A.~Melle--Hern\'andez \thanks{The last two authors were partially
supported by the grant MTM2004-08080-C02-01. Address:
University Complutense de Madrid, Dept. of Algebra,
Madrid, 28040, Spain.
E-mail: iluengo\symbol{'100}mat.ucm.es, amelle\symbol{'100}mat.ucm.es}}
\date{}
\begin{document}
\def\eps{\varepsilon}

\maketitle

\begin{abstract}
We discuss
the notion of a power structure over a ring and
the geometric description of the power structure over the Grothendieck
ring of complex quasi-projective varieties and show some examples of
applications to generating series of classes of configuration spaces
(for example, nested Hilbert schemes of J.~Cheah) and  wreath product
orbifolds.
\end{abstract}

To a pre-$\lambda$ ring there corresponds a so called {\em power
structure}. This means, in particular, that one can give sense to an
expression of the form
$$
\left(1+ a_1 t+a_2+...\right)^m
$$
for $a_i$ and $m$ from the ring $R$. (Generally speaking, on a ring
there are many pre-$\lambda$ structures which correspond to
one and the same
power structure.) A natural pre-$\lambda$ structure on the Grothendieck
ring $K_0(\Var)$ of complex quasi-projective varieties is defined by the
Kapranov zeta-function
$$
\zeta_X(t)=1+[X]t+[S^2X]t2+[S^3X]t3+\ldots\,,
$$
where $S^kX=X^k/S_k$ is the $k$-th symmetric power of the variety $X$.
In \cite{GLM1}, there was given a geometric description of the
corresponding power structure over the Grothendieck ring $K_0(\Var)$. In
some cases this permits to give new (short and somewhat more
transparent) proofs and also certain refinements of formulae for
generating series of classes of moduli spaces in the ring $K_0(\Var)$
and/or of their invariants: the Euler characteristic and the
Hodge--Deligne polynomial. An application of this sort (for the
generating series of classes of Hilbert schemes
of $0$-dimensional subschemes of a smooth quasi-projective variety) was
described in \cite{GLM2}.

The aim of this paper is to describe the concept of a power structure
(in a somewhat more general context introduced in \cite{GLM3}) and to
show its applications to proofs and also some
improvements of
results by J.~Cheah in \cite{Cheah2} about nested Hilbert schemes, by
W.P.~Lin and Zh.~Qin in \cite{LQ1} about moduli spaces of
$1$-dimensional subschemes. Finally
we rewrite some results of W.~Wang, J.~Zhou in \cite{wang} and \cite{wz}
on generating series of orbifold generalized Euler characteristicic of
wreath product orbifolds in terms of the power structure.

\section{Power structures}\label{sec1}
\begin{definition}
A pre-$\lambda$ structure on a ring $R$ is given by a series
$\lambda_a(t)\in 1+t\cdot R[[t]]$ defined for each $a\in R$ so that
\begin{enumerate}
\item $\lambda_a(t)=1+at\, \mbox{ mod }t^2.$
\item $\lambda_{a+b}(t)=\lambda_a(t)\lambda_b(t)$ for $a,b\in R$.
\end{enumerate}
\end{definition}

\begin{example} One has the following important examples of pre-$\lambda$ structures.
\begin{enumerate}
\item $R$ is the ring $\Z$ of integers, $\lambda_k(t)=(1-t)^{-k}.$
\item $R=\Z$ and $\lambda_k(t)=(1+t)^{k}.$
\item $R=\Z[u_1, \ldots, u_r]$ (the ring of polynomials in $r$ variables
$u_1$, \dots, $u_r$), for a polynomial $P=P(\uu)=\sum p_{\kk}\,\uu^{\,\kk}$,
${\kk}\in\Z_{\ge 0}^r$ and $p_{\kk}\in\Z$,
$$
\lambda_{P}(t)=
\prod\limits_{\kk\in\Z_{\ge0}^r} (1-\uu^{\,\kk} t)^{-p_\kk}\,,
$$
where
$\uu=( u_1,\ldots,u_r)$, $\kk=(k_1,\ldots,k_r)$, $\uu^{\,\kk} =u_1^{k_1}\cdot\ldots\cdot u_r^{k_r}$, (see \cite{GLM2}).
\item (A more geometric example.) Let $R$ be the $K$-functor $K(X)$ of the space $X$, i.e. the Grothendieck ring of (say, real or complex) vector bundles over $X$. For a vector bundle $E$ over $X$, let $\Lambda^k E$ be the $k$-th exterior power of the bundle $E$. The series
$$
\lambda_{E}(t)=1+[E]t+[\Lambda^2 E]t^2+[\Lambda^3 E]t^3+\ldots
$$
defines a pre-$\lambda$ structure on the ring $K(X)$.
\end{enumerate}
\end{example}

To a pre-$\lambda$ structure on a ring $R$ one can associate a \emph{power structure} over $R$: a notion introduced in \cite{GLM1}.

\begin{definition}
A {\em power structure} over a (semi)ring $R$ with a unit is a map
$\left(1+t\cdot R[[t]]\right)\times {R} \to 1+t\cdot R[[t]]$:
$(A(t),m)\mapsto \left(A(t)\right)^{m}$,
which possesses the following properties:
\begin{enumerate}
\item[1)] $\left(A(t)\right)^0=1$,
\item[2)] $\left(A(t)\right)^1=A(t)$,
\item[3)] $\left(A(t)\cdot B(t)\right)^{m}=\left(A(t)\right)^{m}\cdot
\left(B(t)\right)^{m}$,
\item[4)] $\left(A(t)\right)^{m+n}=\left(A(t)\right)^{m}\cdot
\left(A(t)\right)^{n}$,
\item[5)] $\left(A(t)\right)^{mn}=\left(\left(A(t)\right)^{n}\right)^{m}$,
\item[6)] $(1+t)^m=1+mt+$ terms of higher degree,
\item[7)] $\left(A(t^k)\right)^m =
\left(A(t)\right)^m\raisebox{-0.5ex}{$\vert$}{}_{t\mapsto t^k}$.
\end{enumerate}
\end{definition}

\begin{remark}
For a ring property 1) follows from the other ones. It is necessary to keep it only for a semiring.
\end{remark}

\begin{definition}
A power structure is {\em finitely determined} if for each $M>0$ there
exists a $N>0$ such that for any series $A(t)$ the $M$-jet of the series $\left(A(t)\right)^m$
(i.e., $\left(A(t)\right)^m\,\mbox{mod }t^{M+1}$) is determined by the $N$-jet
of the series $A(t)$.
\end{definition}

\begin{proposition}\label{prop1}
To define a finitely determined power structure over a ring $R$ it is
sufficient to define the series $(A_0(t))^{m}$ for any fixed series $A_0(t)$ of the form $1+t+$ terms of higher degree, \ and for each $m\in R$, so that:
\begin{enumerate}
\item[1)] $(A_0(t))^m=1+mt+$ terms of higher degree;
\item[2)] $(A_0(t))^{m+n}=(A_0(t))^{m}(A_0(t))^{n}$.
\end{enumerate}
\end{proposition}

\begin{proof}
By properties 6 and 7, each series $A(t)\in 1+t\cdot R[[t]]$ can be written in a
unique way as a product of the form
$\prod\limits_{i=1}^\infty (A_0(t^i))^{b_i}$, with $b_i\in R$. Then by
properties 3 and 7 (and the finite determinacy of the power structure)
\begin{equation}\label{eq1}
\left(A(t)\right)^m=\prod\limits_{i=1}^\infty (A_0(t^i))^{b_i m}.
\end{equation}
In the other direction, one can easily see that the power structure
defined by the
equation~(\ref{eq1})
possesses properties 1)--7).
\end{proof}

Proposition~\ref{prop1} means that a pre-$\lambda$ structure on the ring
$R$ defines a finitely determined power structure over $R$. In the other direction, there are many pre-$\lambda$ structures on the ring $R$ which give one and the same power structure: those defined by the series $(A_0(t))^m$ for any fixed series $A_0(t)$ of the form $1+t+\ldots$ terms of higher degree.
In what follows we shall prefer to use the series $A_0(t)=(1-t)^{-1}=1+t+t^2+\ldots\in R[[t]]$.

\medskip
Let $R[[\tt]]=R[[t_1,\ldots,t_r]]$ be the ring of series in $r$ variables
$t_1$, \ldots, $t_r$ with coefficients from the ring $R$ and let $\gm$ be the
ideal $\langle t_1,\ldots,t_r \rangle$. A power structure over the ring
$R$ in a natural way permits to give sense to expressions of the form
$(A(\tt))^m$, where $A(\tt)\in 1+\gm R[[\tt]]$. Namely, the series
$A(\tt)$ can be in a unique way represented in the form
$$
A(\tt)=\prod\limits_{\kk\in\Z_{\geq 0}^r\setminus \{0\}} (1-\tt^{\,\kk})^{-b_{\kk}}
$$
($\tt^{\kk}=t_1^{k_1}\ldots t_r^{k_r}$).
Then
$$
(A(\tt))^m=\prod\limits_{\kk\in\Z_{\geq 0}^r\setminus \{0\}} (1-\tt^{\,\kk})^{-b_{\kk}m}\,.
$$

Let $R_1$ and $R_2$ be rings with power structures over them. A ring
homomorphism $\varphi:R_1\to R_2$ induces the natural homomorphism
$R_1[[\tt]]\to R_2[[\tt]]$ (also denoted $\varphi$) by
$\varphi\left(\sum a_{\,\ii}\,\tt^{\,\ii}\right)=\sum\varphi(a_{\,\ii})\,\tt^{\,\ii}$.
One has:

\begin{proposition}\label{prop2}
If a ring homomorphism $\varphi:R_1\to R_2$ is such that
$(1-t)^{-\varphi(m)}=\varphi\left((1-t)^{-m}\right)$ for any $m\in R$, then
$\varphi\left(\left(A(\tt)\right)^m\right)=
\left(\varphi\left(A(\tt)\right)\right)^{\varphi(m)}$ for $A(\tt)\in 1+\gm R[[\tt]]$, $m\in R$.
\end{proposition}

\begin{definition}
The {\em Grothendieck ring} $K_0(\Var)$ {\em of complex quasi-projective
varieties} is the abelian group generated by classes $[X]$ of all
quasi-projective varieties $X$ modulo the relations:
\begin{enumerate}
\item[1)] if varieties $X$ and $Y$ are isomorphic, then $[X]=[Y]$;
\item[2)] if $Y$ is a Zariski closed subvariety of $X$, then $[X]=[Y]+[X\setminus Y]$.
\end{enumerate}
The multiplication in $K_0(\Var)$ is defined by the Cartesian product of varieties.
\end{definition}

\begin{remark}
One can also consider the concept of the Grothendieck semiring $S_0(\Var)$ of complex quasi-projective varieties substituting the word ``group" above by the word ``semigroup". Elements of the semiring $S_0(\Var)$ have somewhat more geometric sense: they are represented by ``genuine" quasi-projective varieties (not by virtual ones).
\end{remark}

The class $[\AA^1_\C]\in K_0(\Var)$ of the complex affine line is denoted by $\LLL$. In a number of cases it is reasonable (or rather necessary) to consider the localization $K_0(\Var)[\LLL^{-1}]$ of the Grothendieck ring $K_0(\Var)$ by the class $\LLL$.

For a complex quasi-projective variety $X$, let $S^kX=X^k/S_k$ be the $k$-th symmetric power of the space $X$ (here $S_k$ is the group of permutations on $n$ elements; $S^kX$ is a quasi-projective variety as well).

\begin{definition}
The {\em Kapranov zeta function} of a quasi-projective variety $X$ is the series
$$
\zeta_{X}(t)=1+[X]\cdot t+[S^2X]\cdot t^2+[S^3X]\cdot t^3+\ldots \in K_0(\Var)[[t]]
$$
(\cite{Kap}).
\end{definition}

One can see that
\begin{equation}\label{eq2}
\zeta_{X+Y}(t)=\zeta_{X}(t)\cdot\zeta_{Y}(t)\,.
\end{equation}
This follows from the relation $S^k(X\amalg Y)=\coprod\limits_{i=0}^k{S^iX}\times S^{k-i}Y$.
Also one has
$$
\zeta_{\LLL^n}(t)=\frac{1}{1-\LLL^n t}\,.
$$
As an example this implies that
$$
\zeta_{\CP^n}(t)=\prod\limits_{i=0}^n\frac{1}{1-\LLL^i t}\,.
$$

Equation~(\ref{eq2}) means that the series $\zeta_{X}(t)$ defines a pre-$\lambda$ structure on the
Grothendieck ring $K_0(\Var)$. The geometric description of the corresponding power structures over the ring $K_0(\Var)$ was given in \cite{GLM1}.
We shall formulate it here in the form adapted for series in $r$ variables (\cite{GLM3}).

Let ${A}_{\,\nn}$, $\nn=(n_1,\ldots,n_r)\in \Z^r_{\geq 0}\setminus\{0\}$, and $M$ be
quasi-projective varieties and
$A(\tt)=1+\sum\limits_{{\nn}\in \Z^r_{\geq 0}\setminus\{0\}} [{A}_{\,\nn}]\, {\tt}^{\,\nn}\in
 K_0(\Var)[[\tt]].$
Let $\mathfrak{A}$ be the disjoint union $\coprod\limits_{{\kk}\in \Z^r_{\geq 0}\setminus\{0\}} {A}_{\,\kk}$, and let ${\kk}:\mathfrak{A}\to \Z^r_{\geq 0}$ be the tautological map on it: it sends points of ${A}_{\,\kk}$ to ${\kk}\in \Z^r_{\geq 0}$.

\medskip
{\bf Geometric description of the power structure over the ring $K_0(\Var)$.}
The coefficient at ${\tt}^{\,\nn}$ in the series ${A}(\tt)^{[M]}$ is represented by the configuration space of pairs $(K,\varphi)$, where $K$ is a finite subset of the variety $M$ and $\varphi$ is a map from $K$ to $ \mathfrak{A}$
such that $\sum\limits_{x\in K}{\kk}(\varphi(x))={\nn}$.
To describe such a configuration space as a quasi-projective variety one can write it as
\begin{equation}\label{geom}
\sum_{{\mathbf k}:\,\sum {\ii}\,k_{\,\ii}={\nn}}
\left[
\left(
(\prod_{\ii} M^{k_{\,\ii}})
\setminus\Delta
\right)
\times\prod_{\ii} A_{\ii}^{k_{\,\ii}}/\prod_{\ii} S_{k_{\,\ii}}
\right],
\end{equation}
where ${\mathbf k}=\{k_{\,\ii}:{\ii}\in \Z^r_{\geq 0}\setminus\{0\}, k_{\,\ii}\in\Z\}$ and $\Delta$ is the ``large diagonal" in $M^{\Sigma k_{\,\ii}}$ which consists of $(\sum k_{\,\ii})$-tuples of points of $M$ with at least two coinciding ones; the permutation group $S_{k_{\,\ii}}$ acts by permuting corresponding $k_{\,\ii}$ factors in
$\prod\limits_s M^{k_{\,\ii}}\supset (\prod\limits_{\ii} M^{k_{\,\ii}})\setminus\Delta$ and the spaces $A_{\,\ii}$ simultaneously (the connection between this formula and the description above is clear).

\medskip
One can shows (see \cite{GLM1}) that the described operation really
gives a power structure over $K_0(\Var)$, i.e. it satisfies conditions 1) -- 7)
of the definition. The fact that this structure corresponds to the
Kapranov zeta function follows from the equation
\begin{equation}\label{eq3}
(1+t+t^2+\ldots)^{[M]}=1+[M]\cdot t+[S^2M]\cdot t^2+[S^3M]\cdot t^3+\ldots
\end{equation}
Indeed, since there is only one map from $M$ to a point (the coefficients
in the series $1+t+t^2+\ldots$), the coefficient at $t^n$ in the LHS of
equation (\ref{eq3}) is represented by the space a point of which is
a finite set of points of the variety $M$ with positive multiplicities such that the
sum of these multiplicities is equal to $n$. This is just the definition
of the $n$-th symmetric power of the variety $M$.

It is also useful to describe the binomial $(1+t)^{[M]}$. The coefficient
at $t^n$ in it is represented by the space a point of which is a finite
subset of $M$ with $n$ elements, i.e. the configuration space
$(M^n\setminus \Delta)/S_n$ of unordered $n$-tuples of distinct points
of $M$.

It seems that the power structure can be used to prove some combinatorial identities.
For instance, applying formula (\ref{geom}) to a finite set $M$ with $m$ elements one
gets a formula for the power of a series:
$$\left(1+\sum\limits_{{\nn}\in \Z^r_{\geq 0}\setminus\{0\}} {a}_{\,\nn}\, {\tt}^{\,\nn}\right)^m=
1+\sum\limits_{{\nn}\in \Z^r_{\geq 0}\setminus\{0\}} \left(\sum_{{\mathbf k}:\,\sum {\ii}k_{\ii}={\nn}}
\frac{m!}{(m-\Sigma k_{\ii})!\prod\limits_{\ii} k_{\ii}!}
\prod_{\ii} a_{\ii}^{\,k_{\ii}}\right) \tt^{\,\nn}.
$$

\bigskip
There are two natural homomorphism from the Grothendieck ring  $K_0(\Var)$ to the ring
$\Z$ of integers and to the ring $\Z[u,v]$ of polynomials in two variables: the Euler characteristic (with compact support) $\chi:K_0(\Var)\to \Z$ and the Hodge--Deligne polynomial
$e:K_0(\Var)\to \Z[u,v]$:
 $e(X)(u,v)=
\sum e^{p,q}(X) u^p v^q$.

The formula of I.G.~Macdonald \cite{mac}:
$$
\chi(1+[X]t+[S^2X]t^2+[S^3X]t^3+\ldots)=(1-t)^{\chi(X)}
$$
and the corresponding formula for the Hodge--Deligne polynomial (see \cite{Bur}, \cite[Proposition 1.2]{Cheah1}):
$$
e(1+[X]t+[S^2X]t^2+\ldots)(u,v)=\prod_{p,q}\left(\frac{1}{1-u^p v^q t}\right)^{e^{p,q}(X)}
$$
implies that these homomorphisms respect the power structures on these rings described above
(see Example 3 and Proposition \ref{prop2} or \cite{GLM2}).
Therefore a relation between series from $K_0(\Var)[[\tt]]$ written in terms of the power structure yields  the corresponding relations between the Euler characteristics and the Hodge--Deligne polynomials of these series.

\begin{remark}  It is also possible to define the power structure and to describe it in the \emph{relative setting}, i.e. over the Grothendieck ring $K_0(\VarS)$ of complex quasi-projective varieties over a  variety $S$. The ring  $K_0(\VarS)$ is generated by classes of varieties with
maps ("projections") to $S$. In this case the coefficient of the series $({A}(\tt))^{[M]}$
is the configuration space a point of which is a pair $(K,\varphi)$ consisting of a finite
subset $K\subset M$ which is contained in the preimage of one point of $S$ and the map $\varphi$ commutes with the projections to $S$.
\end{remark}

\section{Nested Hilbert schemes of J.~Cheah}\label{sec2}
Let $\text{Hilb}^n_X$, $n\ge 1$, be the Hilbert scheme of zero-dimensional subschemes of length $n$ of a complex quasi-projective variety $X$; for $x\in X$,
let $\text{Hilb}^n_{X,x}$ be the Hilbert scheme of subschemes of $X$ supported
at the point $x$.

In \cite{Cheah2}, J.~Cheah considered nested Hilbert schemes on a smooth $d$-di\-men\-sio\-nal complex quasi-projective variety $X$.  For $\nn=(n_1,\ldots,n_r)\in\Z_{\geq 0}^r$, the \emph{nested Hilbert scheme} $Z^{\,\nn}_X$ of depth $r$ is the scheme which parametrizes collections of the form
$(Z_1,\ldots,Z_r)$, where $Z_i\in \hilb_X^{n_i}$ and $Z_i$ is a subscheme of $Z_j$ for $i<j$.
The scheme $Z^{\,\nn}_{X}$ is non-empty only if $n_1\leq n_2\leq \ldots\leq n_r$; notice that
$Z^{(n)}_{X}=\hilb_X^{n}\cong Z^{(n,\ldots, n)}_{X}$.

For $Y\subset X$,
let $Z^{\,\nn}_{X,Y}$ be the  scheme which parametrizes collections
$(Z_1,\ldots,Z_r)$ from $Z^{\nn}_X$ with $\text{supp}\,Z_i\subset Y$. For $Y=\{x\}$, $x\in X$, we shall use the notation $Z^{\,\nn}_{X,x}$.

For $r\geq 1$, let $\tt=(t_1,\ldots,t_r)$ and
$$
\ZZ_X^{(r)}(\tt):=\sum\limits_{\nn\in Z_{\geq 0}^r}\,[Z^{\,\nn}_X]\,\tt^{\nn},\quad
\ZZ_{X,x}^{(r)}(\tt):=\sum\limits_{\nn\in Z_{\geq 0}^r}\,[Z^{\,\nn}_{X,x}]\,\tt^{\nn},\quad
$$
be the generating series of classes of nested Hilbert schemes $Z^{\nn}_X$ (resp. supported at the point $x$) of depth $r$.

\begin{theorem}\label{theo1}
For a smooth quasi-projective variety $X$ of dimension $d$, the following identity holds in $S_0(\Var)[[\tt]]$ $($and therefore also in $K_0(\Var)[[\tt]]${}$)$:
\begin{equation}\label{eq4}
\ZZ_X^{(r)}(\tt)=\left(\ZZ_{\AA^d,0}^{(r)}(\tt)\right)^{[X]}.
\end{equation}
\end{theorem}

\begin{proof}
For a Zariski closed subset $Y\subset X$, one has $\ZZ_X^{(r)}(\tt)= \ZZ_{X,Y}^{(r)}(\tt) \cdot \ZZ_{X,X\setminus Y}^{(r)}(\tt)$.
Therefore it is sufficient to prove equation $(\ref{eq4})$ for Zariski open subsets $U$ of $X$ which form a covering of $X$ and for their intersections.

One can take $U$ which lies in an affine chart $\AA^N_\C$ and such
that its projection to a  $d$-dimensional  coordinate space (say, generated by the first $d$ coordinates) is everywhere non-degenerate (i.e. is an \`etale morphism).
For any point $x\in U$, this projection identifies $\nn$-nested Hilbert schemes
of $Z^{\,\nn}_{U,x}$ with $Z^{\,\nn}_{\AA^d_\C,x}$.

A nested (zero-dimensional) subscheme of $U$ of type $\nn$
is defined by a finite subset $K\subset U$ with a nested subscheme from $Z^{\,\kk(x)}_{X,x}$
at each point $x\in K$. such that $\sum\limits_{x\in K}\kk(x)=\nn$.
This coincides with the description of the coefficient at $\tt^{\,\nn}$ in the RHS of the equation (\ref{eq4}).
\end{proof}

Similar considerations permit to give a short proof of a somewhat refined version of the main result of \cite{Cheah2}.
Following J.~Cheah, let
\begin{eqnarray}
{\mathfrak F}_X^n&=&\{(x,Z)\in X\times \hilb_X^n:{\ }x\in \text{supp}\,Z\, \},\nonumber\\
{\mathfrak F}_X^{n-1,n}&=&\{(x_1,x_2,Z_1,Z_2)\in X\times X\times Z^{(n-1,n)}_X:{\ } x_i\in \text{supp}\,Z_i,\,i=1,2 \},\nonumber\\
{\mathfrak T}_X^{n}&=&\{(x_1,x_2,Z)\in X\times X\times \hilb_X^n:{\ } x_i\in \text{supp}\, Z,\,i=1,2 \},\nonumber\\
{\mathfrak G}_X^{n}&=&\{(x,Z_1,Z_2)\in X\times Z^{(n-1,n)}_X:{\ } x\in \text{supp}\,Z_2 \,\}.\nonumber
\end{eqnarray}

Let the series ${\mathfrak P}_X(t_0,t_1, t_2, t_3)$ and ${\mathfrak f}_d(t_0,t_1, t_2, t_3)$ from $K_0(\Var)[[t_0, t_1, t_2, t_3]]$ be defined by
\begin{eqnarray}
{\mathfrak P}_X(t_0,t_1, t_2, t_3)&:=&\left[\sum_{n\ge0}[\text{Hilb}^n_{X}]t_0^n\right]+
\left[\sum_{n\ge1}[{\mathfrak F}_X^n]t_0^n\right]t_1\nonumber\\
&+&
\left[\sum_{n\ge1}[{\mathfrak F}_X^n]t_0^n\right]t_2 + \left[\sum_{n\ge1}[{\mathfrak T}_X^{n}]t_0^n\right]t_1t_2\nonumber\\
&+&
\left[\sum_{n\ge1}[Z^{(n-1,n)}_X]t_0^n\right]t_3 + \left[\sum_{n\ge2}[Z^{(1,n-1,n)}_X]t_0^n\right]t_1t_3\nonumber\\
&+&
\left[\sum_{n\ge1}[{\mathfrak G}_X^{n}]t_0^n\right]t_2t_3 + \left[\sum_{n\ge2}[{\mathfrak F}_X^{n-1,n}]t_0^n\right]t_1t_2t_3\,,\nonumber
\end{eqnarray}
\begin{eqnarray}
{\mathfrak f}_d(t_0,t_1, t_2, t_3)&:=&\sum_{k\ge0}[\hilb^k_{\AA^d,0}]t_0^k + \sum_{k\ge1}[Z^{(k-1,k)}_{\AA^d,0}]t_0^kt_3\nonumber\\
&+&
\sum_{k\ge1}[\hilb^k_{\AA^d,0}]t_0^kt_2 + \sum_{k\ge1}[Z^{(k-1,k)}_{\AA^d,0}]t_0^kt_2t_3\nonumber\\
&+&
\sum_{k\ge1}[\hilb^k_{\AA^d,0}]t_0^kt_1 + \sum_{k\ge1}[\hilb^k_{\AA^d,0}]t_0^kt_1t_2\nonumber\\
&+&
\sum_{k\ge2}[Z^{(k-1,k)}_{\AA^d,0}]t_0^kt_1t_3 + \sum_{k\ge2}[Z^{(k-1,k)}_{\AA^d,0}]t_0^kt_1t_2t_3\,.\nonumber
\end{eqnarray}

\begin{theorem}(cf. Main Theorem in \cite{Cheah2})\label{nested}
Let $X$ be a  smooth quasi-projective variety of dimension $d$. Then
\begin{equation}\label{nest}
{\mathfrak P}_X(t_0,t_1, t_2, t_3)=({\mathfrak f}_d(t_0,t_1, t_2, t_3))^{[X]} \quad  \mod (t_1^2,t_2^2,t_3^2)\,.
\end{equation}
\end{theorem}

\begin{proof}
Using the arguments of the proof of Theorem~\ref{theo1} we may suppose that $X$ lies in an affine chart $\AA^N_\C$ and its projection to a $d$-dimensional coordinate space is
nondegenerate. This identifies $\hilb_{X,x}^s$ and $Z_{X,x}^{(s-1,s)}$ with  $\hilb_{\AA^d,0}^s$ and $Z_{\AA^d,0}^{(s-1,s)}$ respectively for each point $x\in X$. To prove equation (\ref{nest}) one has to give an interpretation of the coefficients at the monomials $t_0^n, t_0^nt_1,\ldots, t_0^nt_1t_2t_3$ in the RHS of (\ref{nest}). Let us make this for the coefficients at $t_0^nt_3$ and at $t_0^nt_2t_3$ (other cases are treated in the same way).

The coefficient at $t_0^nt_3$ is represented by the space a point of which is defined by a point $x_0$ of $X$ with a zero-dimensional nested scheme from $Z_{X,x_0}^{(k(x_0)-1,k(x_0))}$ at it plus several other points of $X$ with  zero-dimensional schemes from $\hilb_{X,x}^{k(x)}\cong Z_{X,x}^{(k(x),k(x))}$ at each of them, such that $k(x_0)+\sum k(x)=n$. This is just the definition of a point of the space $Z_{X}^{(n-1,n)}$.

The monomial $t_0^nt_2t_3$ can be obtained either as a product of two monomials of the form $t_0^*t_2$, $t_0^*t_3$ and of several monomials of the form $t_0^*$ or as a product of a monomial of the form $t_0^*t_2t_3$ and of several monomials of the form $t_0^*$. Therefore the coefficient at the monomial $t_0^nt_2t_3$ is represented by the space consisting of two parts.

A point of the first part is defined by a point $x_1$ of $X$ with a scheme from  $\hilb_{X,x_1}^{k(x_1)}\cong Z_{X,x_1}^{(k(x_1),k(x_1))}$ at it, with $k\geq 1$
(i.e. it is not empty: $x_1$ belongs to the support of it), a point $x_2\in X$ with a scheme from $Z_{X,x_2}^{(k(x_2)-1,k(x_2))}$ at it plus several points of $X$ with $0$-dimensional schemes from  $\hilb_{X,x}^{k(x)}$ at each of them  such that $k(x_1)+k(x_2)+\sum k(x)=n$.

A point of the second part is defined by a point $x_1$ of $X$ with the scheme $(z_1,z_2)$ from
$Z_{X,x_1}^{(k(x_1)-1,k(x_1))}$ at it (in this case $z_2$ is not empty: $x_1$ belongs to the support of it) plus several points of $X$ with $0$-dimensional schemes from  $\hilb_{X,x}^{k(x)}$ at each of them  such that $k(x_1)+\sum k(x)=n$.
Therefore a point of the union of these two subspaces can be described by a nested scheme
$(Z_1,Z_2)$ from $Z_{X}^{(n-1,n)}$ plus a point which belongs to $Z_2$. This is just the description of the space ${\mathfrak G}_X^{n}$.
\end{proof}

Applying the Hodge--Deligne homomorphism to (\ref{nest}) one gets the Main Theorem of \cite{Cheah2}.

\begin{example} \label{incidence}
Let $S$ be a smooth quasi-projective surface. Consider the incidence variety $Z^{(n-1,n)}_{S}=\{ (Z_1,Z_2)\in \hilb_S^{n-1}\times \hilb_S^{n}:\,Z_1\subset Z_2\}$. Using the results of J.~Cheah on the cellular decomposition of $Z^{(n-1,n)}_{\AA^2_\C,0}$ (\cite{Cheah3}), one gets the result of L.~G\"ottsche
(\cite[Theorem 5.1]{Got2}):
\begin{equation*}
\sum_{n\ge 1}[Z^{n-1,n}_S]\,t^n= \frac{[S]\cdot t}{1-\LLL t} \left(\prod\limits_{k\geq 1} \frac{1}{1-\LLL^{k-1}t^k}\right)^{[S]}.
\end{equation*}
\end{example}

\section{On moduli spaces of curves and points \-(after W.-P.~Li and Zh.~Qin)}\label{sec3}
In \cite{LQ1}, there were considered certain moduli spaces of $1$-dimensional subschemes in a smooth $d$-dimensional projective complex variety. Let $X$ be a smooth $d$-dimensional projective complex variety with a Zariski locally trivial fibration $\mu:X\to S$ where $S$ is smooth of dimension $d-1$ and fibres are smooth irreducible curves of genus $g$. Let $\beta\in H_2(X,\Z)$ be the class of the fibre.

Let ${\mathfrak I}_n(X,\beta)$ be the moduli space of $1$-dimensional closed subschemes $Z$ of $X$ such that $\chi({\mathcal O}_Z)=n, [Z]=\beta$, where $[Z]$ is the fundamental class of the scheme $Z$  and let $\gM^n:={\mathfrak I}_{(1-g)+n}(X,\beta).$
Let $\gM^n_{X, C_s, x}$ be the moduli space of $1$-dimensional closed subschemes $\Theta$ in $X$ such that  $I_\Theta\subset I_{C_s}$, the support $\text{supp}\left({I_{C_s}}/{I_\Theta}\right)=\{O\}$ and $\dim_O\left({I_{C_s}}/{I_\Theta}\right)=n$. The number $n$ will be called the length of the subscheme $\Theta$. Let $\glocal^n$ have the same meaning: it is the moduli space of $1$-dimensional closed subschemes $\Theta$ in $\C^{d-1}\times \C$ such that  $I_\Theta\subset I_{\{O\}\times \C}$, $\text{supp}\left({I_{\{O\}\times \C}}/{I_\Theta}\right)=\{O\}$ and $\dim_O\left({I_{\{O\}\times \C}}/{I_\Theta}\right)=n$.

\begin{theorem}(cf. Proposition 5.3, Lemma 6.1 and Proposition 6.2  in \cite{LQ1})\label{moduli}
Let $X$ be a smooth $d$-dimensional projective complex variety with a Zariski locally trivial fibration $\mu:X\to S$ where $S$ is smooth of dimension $d-1$ and fibres are smooth irreducible curves of genus $g$. Then
\begin{equation}\label{eqmoduli}
\sum_{n\ge 0}[\gM^n]\,t^n=[S] \left(\sum_{n\ge 0}[\hilb^n_{\C^d,0}]\,t^n \right)^{[X]-[C]}\left(\sum_{n\ge 0}[\glocal^n]t^n\right)^{[C]}.
\end{equation}
\end{theorem}

\begin{proof}
A point of $\gM^n$ can be considered as consisting of a fibre $C_s=\mu^{-1}(s)$ of the bundle $\mu:X\to S$ and of several fixed points, both outside of $C_s$ and on it, with a $0$-dimensional subscheme (i.e. an element of $\hilb_{X,x}^*$) at each of those points which are outside of $C_s$ and a subscheme of ${\mathfrak M}_{X,C_s,x}^*$ at each of those points which lie on $C_s$ such that  the sum of their lengths is equal to $n$. Thus, there is a natural map (projection) from $\gM^n$ to $S$.
Over a point $s\in S$, there are somewhat different objects (subschemes) at points outside of the curve $C_s$ and on this curve.

It is sufficient to prove equation (\ref{eqmoduli}) for preimages of elements of a covering of $S$ by Zariski open subsets and of their intersections. Therefore without any loss of generality we can suppose that $X=S\times C$.
Moreover, let us choose a fixed point $s_0\in S$. A constructible map which sends ${\mathfrak M}^n_{X,C_s}$ to  ${\mathfrak M}^n_{X,C_{s_0}}$ and is an isomorphism of strata can be defined as follows. One takes a $0$-dimensional subscheme which lies on $C_{s_0}$ and puts them to the corresponding points of $C_s$ and vice versa, one takes the elements of ${\mathfrak M}^n_{X,C_{s},x}$ and puts them to the corresponding points  of $C_{s_0}$. Thus in the Grothendieck ring of algebraic varieties one has $[{\mathfrak M}^n]=[S][{\mathfrak M}^n_{X,C_{s_0}}]$.

Therefore to prove (\ref{eqmoduli}) one should show that
\begin{equation}\label{eqmoduli2}
\sum_{n\ge 0}[{\mathfrak M}^n_{X,C_{s_0}}]\,t^n=\left(\sum_{n\ge 0}[\hilb^n_{\C^d,0}]\,t^n \right)^{[X]-[C_{s_0}]}\left(\sum_{n\ge 0}[\glocal^n]\,t^n\right)^{[C_{s_0}]}.
\end{equation}
Just as in the proofs in Section~\ref{sec2} we may suppose that at each point of $X$ the space $\hilb^k_{X,x}$ is identified with the space $\hilb^k_{\C^d,0}$ and at each point of $C_{s_0}\subset X$ the space ${\mathfrak M}^k_{X,C_{s_0}, x}$ is identified with the space ${\mathfrak M}^k_{\C^{d-1}\times\C, \{0\}\times\C, 0}$. The coefficient at $t^n$ in the RHS of equation~(\ref{eqmoduli2}) is represented by the space a point of which is defined by several points of the curve $C_{s_0}\subset X$ with a scheme from ${\mathfrak M}^{k(x)}_{X,C_{s_0}, x}$ at each of them and several points from $X\setminus C_{s_0}$ with a scheme from $\hilb^{k(x)}_{X,x}$ at each of them such that the sum of the lengths $k(x)$ over all the mentioned points is equal to $n$. This is just the description of a point of ${\mathfrak M}^n_{X,C_{s_0}}$.
\end{proof}

\section{Orbifold generalized Euler characteristic and the power structure}\label{sec4}
Here we rewrite some results of \cite{wang} and \cite{wz} in terms of the power structure. For that we need it over a somewhat modified version of the Grothendieck ring $K_0(\Var)$. For a fixed positive integer $m$, consider the ring $K_0(\Var)[\LLL^{1/m}]$. The pre-$\lambda$ structure on (and therefore the corresponding power structure over) the ring $K_0(\Var)$ can be extended to one on $K_0(\Var)[\LLL^{1/m}]$ by the formula
$$
\zeta_{[X]\LLL^{s/m}}(t)=\zeta_{[X]}(\LLL^{s/m}t)\,.
$$
In a similar way the corresponding pre-$\lambda$ structure on the ring $\Z[u_1^{1/m}, \ldots, u_r^{1/m}]$
can be defined by the formula
$$
\lambda_{P}(t)=
\prod\limits_{\kk\in\Z_{\ge0}^r} (1-\uu^\kk t)^{-p_{\kk}}
$$
for a polynomial $P=P(\uu)=\sum\limits_{\kk\in(1/m)\Z_{\ge0}^r} p_{\,\kk}\,\uu^{\,\kk}$.
There are natural homomorphisms ($\chi$ and $e$) from the ring $K_0(\Var)[\LLL^{1/m}]$ to the rings $\Z$ and $\Z[u^{1/m}, v^{1/m}]$ which send the element $\LLL^{1/m}$ to $1$ and $(uv)^{(1/m)}$ respectively. One can easy see that these are homomorphisms of the $\lambda$-rings and therefore they respect the power structures.

\bigskip

Let $X$ be a smooth quasi-projective complex algebraic variety of dimension $d$ with an action of a finite group $G$ of order $m$. For an element $g\in G$, let $X^g$ be the set $\{x\in X: gx=x\}$ of $g$-invariant points of the action. If $h=vgv^{-1}$ in $G$, the element $v$ defines an isomorphism $v:X^g\to X^h$. Let $G_*$ be the set of conjugacy classes of elements of the group $G$. For a conjugacy class $c\in G_*$ choose its representative $g\in G$.  Let $C_G(g)$ be the centralizer of the element $g$ in $G$.  The centralizer $C_G(g)$ acts on the set $X^g$ of fixed points of $g$. Suppose that its action on the set of connected components of $X^g$ has $N_c$ orbits and let $X^g_1$, \dots, $X^g_{N_c}$ be unions of components of each of these orbits. At each point $x\in X_{\alpha_c}^g$, the map $dg$ is an automorphism of the tangent space $T_x X$ which acts as a diagonal matrix $\mbox{diag\,}(\exp(2\pi i \theta_1),\ldots,\exp(2\pi i \theta_d))$, where $0\leq \theta_i<1$, $\theta!
 _i\in (1/m)\Z$. The \emph{shift number} $F_{\alpha_c}^g$ associated to $X_{\alpha_c}^g$ is
$F_{\alpha_c}^g:=\sum_{j=1}^d \theta_j\in \Z/m$ (it was introduced by E.~Zaslow
in \cite{Zas}).

\begin{definition} The {\em orbifold generalized Euler characteristic} $[X,G]$
of the pair $(X,G)$ is
$$
[X,G]:=\sum_{c\in G_*} \sum_{\alpha_c=1}^{N_c} [X_{\alpha_c}^g/C_G(g)]\cdot\LLL^{F_{\alpha_c}^g} \in K_0(\Var)[\LLL^{1/m}].
$$
\end{definition}

Applying the Euler characteristic morphism one gets the notion of orbifold Euler characteristic invented in the study of string theory of orbifolds by L.~Dixon et al. \cite{Dixon}:
$$
\chi (X,G):=\sum_{c\in G_*} \sum_{\alpha_c=1}^{N_c} \chi(X_{\alpha_c}^g/C_G(g))=\sum_{c\in G_*} \chi(X^g/C_G(g))\,.
$$
Applying the Hodge--Deligne polynomial one gets the orbifold $E$-function
introduced by V.~Batyrev in \cite{baty}:
$$
E_{orb}(X,G;u,v):=\sum_{c\in G_*} \sum_{\alpha_c=1}^{N_c} e(X_{\alpha_c}^g/C_G(g))(u,v)\, (uv)^{F_{\alpha_c}^g} \in \Z[u^{1/m},v^{1/m}]\,.
$$

Let $G^n=G\times\ldots\times G$ be the Cartesian power of the group $G$. The symmetric group $S_n$ acts on $G^n$ by permutation of the factors: $s (g_1, \ldots, g_n) = (g_{s^{-1}(1)} , \ldots, g_{s^{-1}(n)})$. The {\em wreath product}
$G_n = G \sim S_n$ is the semidirect product of $G^n$ and $S_n$ defined by the described action. Namely the multiplication in the group $G_n$ is given by the formula $(g, s)(h, t) = (g\cdot s(h), st)$, where $g,\, h \in G^n$,
$s,\, t \in S_n$. The group $G^n$ is a normal subgroup of $ G_n$ via the identification of $g \in G^n$ with $(g,1) \in G_n$.
For a variety $X$ with a $G$-action, there is the corresponding action of the group $G_n$ on the Cartesian power $X^n$ given by the formula
$$
( (g_1, \ldots, g_n), s)(x_1, \ldots, x_n) = (g_1 x_{s^{-1} (1)}, \ldots, g_n x_{s^{-1} (n)})\,,
$$
where $x_1, \ldots, x_n \in X$, $g_1, \ldots, g_n \in G$, $s\in S_n$. One can see that the factor variety $X^n/G_n$ is naturally isomorphic to $(X/G)^n/S_n$. In particular, $[X^n/G_n]=[(X/G)^n/S_n]$ in $ K_0({\Var})$. Therefore
\begin{equation*}
\sum_{n \geq 0} [X^n/G_n] t^n = (1-t)^{-[X/G]}\in K_0({\rm{\Var}})[[t]]\,.
\end{equation*}

\begin{theorem} (cf. \cite{wang},\cite{wz})
Let $X$ be  a smooth quasi-projective complex algebraic variety of dimension $d$ with an action of a finite group $G$ of order $m$. Then
\begin{equation}\label{eq5}
\sum_{n \geq 0} [X^n, G_n] t^n =  \left(\prod_{ r =1}^{\infty}
( 1 -\LLL^{(r-1)d/2}t^r)\right)^{-[X, G]}\,.
%% = \left(\Ex \left(\frac{\LLL^{d/2}t}{1-\LLL^{d/2}T}\right)\right)^{\LLL^{(-d/2)}[X, G]}
\end{equation}
\end{theorem}

\begin{proof}
One can say that essentially the proof is already contained in \cite{wz} where invariants of the $G_n$-action on the space $X^n$ are related to those of the $G$-action on the space $X$ (see also \cite{wang} and \cite{tamanoi}).

Let $a=(g,s)\in G_n$, $g=(g_1,\ldots, g_n)$. Let $z=(i_1,\ldots, i_r)$ be one of the cycles in the permutation $s$. The {\em cycle-product} of the element $a$ corresponding to the cycle $z$ is the product $g_{i_r}g_{i_{r-1}}\ldots g_{i_1}\in G$. The conjugacy class of the cycle-product is well-defined by $g$ and $s$. For $c\in G_*$ and $r\ge 0$, let $m_r(c)$ be the number of $r$-cycles in the permutation $s$ whose cycle-products lie in $c$. Let $\rho(c)$ be the partition which has $m_r(c)$ summands equal to $r$, and let $\rho=(\rho(c))_{c\in G_*}$ be the corresponding partition-valued function on $G_*$. One has
$$
\Vert\rho \Vert : = \sum\limits_{c \in G_*} |\rho (c)| = \sum\limits_{c \in G_*, r\geq 1} r m_r (c) = n\,.
$$
The function $\rho$ or, equivalently, the data $\{m_r(c)\}_{r,c}$ is called the {\em type} of the element $a=(g,s)\in G_n$.
Two elements of the group $G_n$ are conjugate to each other iff they are of the same type.

In \cite{wz} it is shown that:
\begin{enumerate}
\item
For a conjugacy class of elements of the group $G_n$ containing an element $a$ of type $\rho=\{m_r(c)\}_{r\ge1,\,c\in G_*}$ ($\sum_{r,c}rm_r(c)=n$), the subspace $(X^n )^a$ can be naturally identified with $\prod\limits_{c,r} (X^c)^{m_r (c)}$.
The factor space $(X^n )^a /Z_{G_n}(a)$ is naturally isomorphic to $\prod\limits_{c \in G_*, r \geq 1} S^{m_r (c)} \left( X^{c} / Z_G (c) \right)$. The connected components of the space $(X^n )^a /Z_{G_n}(a)$ are numbered by sets of integers $(m_{r,c}(1), \ldots, m_{r,c}(N_c))$ satisfying the relation $\sum\limits_{\alpha_c =1}^{N_c} m_{r,c}(\alpha_c) = m_{r}(c)$.
They are
$$
(X^n)^a_{\{m_{r,c}(\alpha_c)\}} =
\prod\limits_{c \in G_*, r \geq 1} \prod_{\alpha_c =1}^{N_c} S^{m_{r,c}(\alpha_c)}(X^c_{\alpha_c} /Z_G(c))\,.
 $$
\item
The shift for the component $(X^n)^a_{\{m_{r,c}(\alpha_c)\}}$ is equal to
$$
F_{\{m_{r,c}(\alpha_c)\}}
 = \sum_{c \in G_*, r \geq 1} \sum_{\alpha_c =1}^{N_c} m_{r,c}(\alpha_c) \left(
F^c_{\alpha_c}(r-1)d/2\right)\,.
$$
\end{enumerate}

These two facts imply that
\begin{eqnarray}
\sum_{n \geq 0} [X^n, G_n] t^n&=&\sum_{n \geq 0} \left( \sum_{m_r(c)}
\prod_{c ,r}\prod_{\alpha_c=1}^{N_c}[S^{m_{r,c}}(X_{\alpha_c}^g/Z_G(g))]
\LLL^{m_r(c)(F_{\alpha_c}^g+\frac{(r-1)d}{2})} \right) t^n\nonumber\\
\allowbreak
&=&
\sum\limits_{ m_r(c)}\prod\limits_{c,r}
   \left(\prod\limits_{\alpha_c=1}^{N_c}[S^{m_{r,c}(\alpha_c)}\left(X_{\alpha_c}^g/Z_G(g)\right)
\LLL^{m_r(c)(F_{\alpha_c}^g+\frac{(r-1)d}{2})} \right)  t^{r{m_r(c)}}\nonumber\\
\allowbreak
&=&
\prod\limits_{c,r}  \prod\limits_{\alpha_c=1}^{N_c} \left(\sum_{m_{r,c}(\alpha_c)}[S^{m_{r,c}(\alpha_c)}\left(X_{\alpha_c}^g/Z_G(g)\right)]
\LLL^{m_r(c)(F_{\alpha_c}^g+\frac{(r-1)d}{2})}  t^{rm_{r,c}(\alpha_c)}\right)\nonumber\\
\allowbreak
&=&
\prod\limits_{c,r}  \prod\limits_{\alpha_c=1}^{N_c} \left(1-\LLL^{(F_{\alpha_c}^g+\frac{(r-1)d}{2})}  t^r\right)^{-[X_{\alpha_c}^g/Z_G(g)]}\nonumber\\
\allowbreak
&=&
\prod\limits_{c,r}  \prod\limits_{\alpha_c=1}^{N_c} \left(1-\LLL^{\frac{(r-1)d}{2})}  t^r\right)^{-\LLL^{F_{\alpha_c}^g}[X_{\alpha_c}^g/Z_G(g)]}\nonumber\\
\allowbreak
&=&
\prod\limits_{r\geq 1}  \left(1-\LLL^{\frac{(r-1)d}{2})}  t^r\right)^{-[X,G]} =\prod\limits_{r\geq 1}  \left(1-(\LLL^{\frac{d}{2}} t)^r\right)^{-\LLL^{-d/2}[X,G]}\,.\nonumber
\end{eqnarray}
\end{proof}

Taking the Euler characteristic of the both sides of the equation (\ref{eq5}), one gets Theorem~5 of
\cite{wang}:
$$
\sum_{n \geq 0} \chi (X^n, G_n)\, t^n =  \prod_{ r =1}^{\infty}
  ( 1 -t^r)^{- \chi(X, G)}.
$$

Applying the Hodge--Deligne polynomial homomorphism, one gets the main result of \cite{wz}:
$$
\sum_{n=1}^{\infty} e(X^n,G_n; u,v)\, t^n =
\prod_{r=1}^{\infty}\prod_{p,q}\left(\frac{1} {(1 - u^pv^q
t^r(uv)^{(r-1)d/2})}\right)^{e^{p,q}_{(X,G)}}\,.
$$

\end{document}